\documentclass[12pt]{amsart}

\newcommand{\rk}{\operatorname{rank}}
\newcommand{\en}{\operatorname{End}}

\newtheorem{Theorem}{Theorem}[section]
\newtheorem{Proposition}[Theorem]{Proposition}
\newtheorem{Definition}[Theorem]{Definition}
\newtheorem{Lemma}[Theorem]{Lemma}\newtheorem{Example}[Theorem]{Example}
\newtheorem{Corollary}[Theorem]{Corollary}

\newcommand{\el}{\mbox{$\mathcal L$}}
\newcommand{\dee}{\mbox{$\mathcal D$}}
\newcommand{\sqs}{\mbox{${\mathcal S}(S)$}}

\newcommand{\jay}{\mbox{$\mathcal J$}}
\newcommand{\ar}{\mbox{$\mathcal R$}}
\newcommand{\eh}{\mbox{$\mathcal H$}}
\newcommand{\leqel}{\mbox{$\leq _{\mathcal L }$}}
\newcommand{\leqar}{\mbox{$\leq _{\mathcal R}$}}
\newcommand{\leqels}{\mbox{$\leq _{{\mathcal L}^{\ast}}$}}
\newcommand{\leqars}{\mbox{$\leq _{{\mathcal R}^{\ast}}$}}
\newcommand{\leqjays}{\mbox{$\leq _{{\mathcal J}^{\ast}}$}}
\newcommand{\leqjay}{\mbox{$\leq _{{\mathcal J}}$}}

\newcommand{\els}{\mbox{${\mathcal L}^{\ast}$}}
\newcommand{\jays}{\mbox{${\mathcal J}^{\ast}$}}

\newcommand{\ars}{\mbox{${\mathcal R}^{\ast}$}}
\newcommand{\ehs}{\mbox{${\mathcal H}^{\ast}$}}
\newcommand{\dees}{\mbox{${\mathcal D}^{\ast}$}}
\renewcommand{\a}{\mbox{$\alpha$}}

\newcommand{\leql}{\mbox{$\leq _{\ell }$}}
\newcommand{\leqr}{\mbox{$\leq _r$}}
\newcommand{\leqj}{\mbox{$\leq _j$}}
\newcommand{\s}{\mbox{$^{\sharp}$}}

\begin{document}

\title[Semigroups of quotients]{Semigroups of left quotients - the layered
approach}

\keywords{Group inverse, order, partial order, quotients
\\
\hspace*{4mm}2000 {\em Mathematics Subject Classification} 20 M 07, 20 M 30}

\author{Victoria Gould}
\address{Department of Mathematics\\University of York\\Heslington\\York
YO10 5DD\\UK}
\email{varg1@york.ac.uk}
\date{\today}

\begin{abstract}

A subsemigroup $S$ of a semigroup $Q$ is a {\em left order}
in $Q$ and $Q$ is a {\em semigroup of left quotients} of
$S$ if every element of $Q$ can be expressed as $a\s b$ where
$a,b\in S$ and if, in addition, every element of $S$ that
is {\em square cancellable} lies in a subgroup of $Q$.
Here  $a\s$ denotes the inverse of $a$ in a subgroup
of $Q$. We say
that a left order $S$ is {\em straight}
in $Q$ if in the above definition
we can insist that $a\,\ar\, b$ in $Q$. A 
complete characterisation
of straight left orders in terms of {\em embeddable $*$-pairs}
is available. In this paper we adopt a different approach,
based on {\em partial order decompositions} of semigroups. 
Such decompositions include semilattice decompositions
and decompositions of a semigroup into principal factors or
principal $*$-factors. We determine when a semigroup that
can be decomposed into straight left orders is itself a
straight left order. This technique gives
a unified approach to obtaining many of
the early results on characterisations of straight left
orders.
\end{abstract}
\maketitle

\section{Introduction}

We assume the reader has a familiarity with Green's relations
$\el,\ar,\eh,\dee$ and $\jay$ and
the preorders $\leqel, \leqar$ and $\leqjay$,
 as may be found in \cite{howie}, and also with
the `starred' generalisations $\els,\leqels$ etc.
 The best reference
for these is \cite{fountainabundant}.
For further detail concerning left orders we refer the reader
to \cite{giiigiv}.

The notion of a left order in a semigroup was introduced
by Fountain and Petrich in \cite{fp} and has been
widely used since. Fountain and Petrich describe those
semigroups that are left orders in completely 0-simple
semigroups.
 We remark here that if
$S$ is a left order in $Q$, $\eh$ is a congruence on $Q$ and
$Q$ is regular, then $S$ is perforce straight \cite{vienna}.
Subsequently, a number of papers have appeared characterising
semigroups that are (straight) left orders in semigroups
in various classes, for example \cite{clifford, compreg}. 
More recently, 
 Theorem 4.4 of \cite{giiigiv} gives 
a characterisation of straight left orders
in arbitrary semigroups; 
we recap briefly below the
approach of \cite{giiigiv}.

An ordered pair $(\leql,\leqr)$ of preorders on a semigroup $S$ is a
$*$-{\em pair} if $\leq_l$ is right compatible with
multiplication, $\leq_r$ is left compatible, $\leq_l\,\subseteq\,\leqels$
and $\leq_r\,\subseteq\,\leqars$.  Clearly 
$(\leqels,\leqars)$ is
a $*$-pair for any semigroup $S$.

\begin{Lemma}\label{induced}
Let $S$ be a subsemigroup of $Q$. Then
\[(\leq_{\mathcal{L}^Q}\cap \,(S\times S),
\leq_{\mathcal{R}^Q}\cap \,(S\times S))\]
is a $*$-pair for $S$, the $*$-{\em pair for $S$ induced by $Q$.}
\end{Lemma}

In particular the above applies when $S$ is a left order in $Q$.

For a $*$-pair $(\leql,\leqr)$, 
we denote by $\el'$ and $\ar'$ the
equivalence
relations associated with $\leq_l$ and $\leq_r$ respectively. Notice
that $\el'$ is a right congruence contained in $\els$ and $\ar'$
is a left congruence contained in $\ars$; $\eh'$ denotes
the intersection of $\el'$ and $\ar'$ and we put
\[\mathcal{G}(S)=\{ a\in S:a\,\eh'\, a^2\} .\]
Notice that $\mathcal{G}(S)\subseteq\sqs$,
where
\[\sqs=\{ a\in S:a\,\ehs\, a^2\}\]
is the set of {\em square cancellable} elements of $S$.
We put $\dee'=\el'\vee\ar'$ and define the relation
$\leqj$ by the rule that for any $a,b\in S$,
\[a\,\leqj\, b\mbox{ if and only if }a\,\dee'\, ubv\mbox{ for some
}u,v\in S.\]
Finally we say that for any $a,b\in S$
\[a\,\jay'\, b\mbox{ if and only if }a\,\leqj\,b\,\leqj\, a,\]
so that if $\leqj$ is a preorder, then $\jay'$ is the
associated equivalence relation.

In the following definition, the slightly eccentric notation is that
inherited
from \cite{giiigiv}.

\begin{Definition}\cite{giiigiv} Let
$(\leql,\leqr)$ be a $*$-pair
for  $S$. Then $(\leql,\leqr)$ is
an {\em embeddable $*$-pair} if the following conditions and the duals
(Eii)(r), (Ev)(r), (Evi)(r) and (Evii)(r)  of
(Eii)(l), (Ev)(l), (Evi)(l) and (Evii)(l) respectively hold.

(Ei) $\el'\circ\ar'=\ar'\circ\el'$.

(Eii)(l) For all $b,c\in S$, $b\leq_l c$ if and only if $b\,\el'\, dc$
for some $d\in S$.

(Eiii)  Every $\el'$-class and every $\ar'$-class contains an element
from $\mathcal{G}(S)$.

(Ev)(l) For all $a\in\mathcal{G}(S)$ and $b\in S$, if $b\leq_l a$, then
$ba\,\ar'\, b$.

(Evi)(l) For all $a\in\mathcal{G}(S)$ and $b,c\in S$, if $b,c\leq_l a$
and $ba=ca$, then $b=c$.

(Evii)(l) For all $a\in\mathcal{G}(S)$ and $b,c\in S$, if $b,c\leq_l a$
and $ba\,\el'\, ca$, then $b\,\el'\, c$.
\end{Definition}

 Embeddable $*$-pairs
are  crucial in determining the existence of semigroups of
straight left quotients, as the next result shows.
We recall that a subset $T$ of a semigroup $S$ is
{\em right reversible}, or {\em satisfies the left Ore
condition}, if for any $a,b\in T$ there exist $c,d\in T$ such
that $ca=db$.

\begin{Theorem}\label{biggie}
\cite{giiigiv} Let $S$ be a semigroup having an embeddable $*$-pair
$(\leql,\leqr)$.
Then
$S$ is a straight left order in a semigroup $Q$ inducing
$(\leql,\leqr)$ if and only if $(\leql,\leqr)$
satisfies the following conditions:

(Gi) $\sqs =\mathcal{G}(S)$;\\
and

(Gii) If $a\in \sqs$ then $H'_a$ is right reversible.
\end{Theorem}

In some sense Theorem~\ref{biggie} is the final answer in the
characterisation of straight left orders. However, earlier
descriptions of straight left orders in semigroups from certain
classes, for example in Clifford semigroups, concentrated
on decomposing a semigroup $S$ into subsemigroups
 that
are straight left orders, and then glueing together the
semigroups of left quotients of the subsemigroups of $S$ to
form a semigroup of left quotients of
the entire semigroup $S$. For example,
Theorem 3.1 of \cite{clifford} states that $S$ is a left
order in a semilattice $Y$ of groups $G_{\alpha},
\alpha\in Y$, if and only if $S$ is a semilattice $Y$
of right reversible, cancellative semigroups $S_{\alpha},\a\in Y$.
The latter are known by the theorem of Ore and Dubreil to be
left orders in groups \cite{cp}. A similar approach is
taken to straight left orders in completely simple semigroups
in \cite{compreg}.
The approach of \cite{clifford,compreg}
and other early papers was thus to determine the global
structure of a semigroup from the local structure. Other
results, for example Theorem 5.12 of \cite{basisii},
do not use this approach explicitly, but nevertheless,
the semigroups in question have a relative ideal series
for which the Rees quotients are left orders in completely
0-simple semigroups.

The aim of this paper is to investigate the relation between
a semigroup being a straight left order, and its decompositions
into subsemigroups or Rees quotients
 that are straight left orders. To this end
we introduce the idea of a {\em partial order decomposition}
of a semigroup, that allows us to consider at the same
time semilattice decompositions and decompositions into
Rees quotients.  This technique enables us to extend the
main theorem of \cite{compreg} to arbitrary semilattices
of orders, find an alternative proof of Theorem 5.12
of \cite{basisii}, and a new characterisation of straight
left orders in completely semisimple semigroups.

The structure of the paper is as follows. In  Section ~\ref{preliminaries}
we prove a number of minor technical results needed for the
remainder of the paper. In Section ~\ref{po} we introduce
the notion of a decomposition of a semigroup $S$ as a partial
order $P$ of subsets $S_{\alpha},\a\in P$. Associated
with each $\alpha\in P$ is a Rees quotient 
$\overline{S_{\alpha}}$.

 Section ~\ref{slicing} looks
at the situation where $S$ is a straight left order in $Q$ and
$Q$ is a partial order $P$ of subsets $Q_{\alpha},\a\in P$.
From the definition of partial order of
subsets we have that $S$ is a partial order
$P$  of subsets $S_{\alpha}=Q_{\alpha}\cap S$.
We consider under what circumstances the associated
Rees quotients $\overline{S_{\alpha}}$
are  straight left orders.
We investigate the relation between $P$ being 
associated with the $\leqjay$ order on $Q$ and
$P$ being associated with the $\leqjays$ order on $S$.
Section ~\ref{layering}
 concentrates on the converse to the scenario
in Section ~\ref{slicing}, namely, given a semigroup $S$ that is a partial
order $P$ of subsets $S_{\alpha},\a\in P$,
such each $\overline{S_{\alpha}}$ is a 
straight left order, when is $S$ a straight left order?
Theorem ~\ref{main} gives necessary and sufficient conditions
for this to be the case.

The remaining sections of the paper are specialisations
of Theorem ~\ref{main}. In Section ~\ref{semilattices}
 we consider semilattices
of straight left orders. In Section ~\ref{fs} we look at straight left
orders that are {\em fully stratified}, that is, the $*$-pair
induced by the semigroup of left quotients is $(\leqels,\leqars)$.
In the case where $S$ is {\em idempotent connected abundant},
we obtain an appropriately clean version of the results of
Sections ~\ref{slicing} and ~\ref{layering}. This enables us to show that the semigroup
of endomorphisms of finite rank of a stable basis algebra
are a straight left order. The final section gives a new
characterisation of straight left orders in completely semisimple
semigroups.

\section{Preliminaries}\label{preliminaries}

We begin with a technical result on zeroes.

\begin{Lemma}\label{zeroes} Let $S$ be a straight left order
in $Q$. Then

(A)  $S$ has a zero if and only if
$Q$ has a zero;

(B) if $S$ has a zero, then $S$ has no zero divisors if and only
if $Q$ has no zero divisors.

\end{Lemma}
\begin{proof} (A) If $Q$ has a zero $0$, then as $0\,\eh\, x$ for some
$x\in S$, we must have that $0=x\in S$. The converse
follows from Proposition 2.8 of \cite{hii}.

(B) Suppose that $S$ and $Q$ have a zero $0$ and 
that $S$ has no zero divisors. Let $p,q\in Q\setminus\{ 0\}$
with $p\,\eh\, a,q\,\eh\, b$ and $a,b\in S\setminus\{ 0\}$.
Then
\[pq\,\ar\, pb\,\el\, ab\]
so that as $ab\neq 0$ we have $pq\neq 0$. Hence
$Q$ has no zero divisors. The converse is clear.
\end{proof}

We recall that for a semigroup $S$, $S^0$ denotes $S$
with a zero adjoined. The following observation will be
useful.

\begin{Lemma}\label{zeroesagain} Let $S$ be a semigroup.
Then $S$ is a straight left order in $Q$ if and only if $S^0$
is a straight left order in $Q^0$.
\end{Lemma}

Next we recall from \cite{basisii} that a semigroup
$S$ is {\em abundant} if every $\els$-class and every
$\ars$-class of $S$ contains an idempotent.
In a non-regular abundant semigroup $S$, the relationship
between the idempotents in $R^*_a$ and those in
$L^*_a$ for a given $a\in S$ need not be as strong
as it is in a regular semigroup. In view of
this, El-Qallali and Fountain \cite{elqallalifountain} introduced
the notion of an  {\em idempotent connected}, or IC,
abundant semigroup. An abundant semigroup is IC
if for each $a\in S$ and each idempotent $e\,\leqels\, a$
($e\,\leqars\, a$), there is an element $b\in S$
with
$ae=ba$ ($ea=ab$) \cite{armstrong}.

\begin{Lemma}\label{ic} Let $S$ be an IC abundant semigroup. Then for
any $a,b\in S$
\[a\,\leqels\, b\mbox{ if and only if }a\,\els\, cb\mbox{ some }c\in
S,\]
and dually for $\leqars$.
\end{Lemma}
\begin{proof} If $a\,\els\, cb$, then it follows from the definitions
that $a\,\leqels\, b$.

Conversely, suppose that $a\,\leqels \, b$. Now
$a\,\els\, e$ and $b\,\els\, f$ for some $e,f\in E(S)$.
Hence $e\,\leqels\, f$, whence $e=ef$ and
\[a\,\els\, e\,\el\, fe\,\els\, be.\]
But $S$ is IC abundant, so that $be=cb$ for some $c\in S$.
\end{proof}

Finally in this section we consider the restrictions
of $\dee, \jay$ and $\leqjay$ on a semigroup of left quotients
to a left order.

\begin{Lemma}\label{jdonS} Let $S$ be a straight left order
in $Q$ inducing $(\leql,\leqr)$. Then
\[\dee'=\el'\circ\ar'=\ar'\circ\el',\]
consequently,
\[\dee'=\dee\cap (S\times S).\]
Further,
\[ \leqj=\leqjay\cap(S\times S) \mbox{ and }
\jay'=\jay\cap(S\times S).\]
\end{Lemma}
\begin{proof} If $a,b\in S$ and $a\, \el'\, c\,\ar'\, b$
for some $c\in S$, then as $Q$ induces $(\leql,\leqr)$
we have that $a\,\dee\, b$ in $Q$. Now $S$ intersects every
$\eh$-class of $Q$, whence $a\,\ar'\, d\,\el'\, b$
for some $d\in S$. Together with the dual argument we
have that $\el'$ and $\ar'$ commute, so that
\[\dee'=\el'\circ\ar'=\ar'\circ\el'\]
as required. Clearly then $\dee'=\dee\cap(S\times S)$.

Suppose now that $a,b\in S$ and $a\,\leqj\, b$. By
definition, we have $a\,\dee'\, ubv$ for some $u,v\in S$.
From the previous paragraph and the fact that $\dee\subseteq \jay$,
\[a\,\jay\, ubv\,\leqjay\, b\]
in $Q$.

Conversely, let $a,b\in S$ and suppose that $a\,\leqjay\, b$ in $Q$.
Then $a=pbq$ for some $p,q\in Q$. Choose $u,v\in S$ with
$p\,\eh\, u$ and $q\,\eh\, v$ so that
\[a=pbq\,\ar\, pbv\,\el\, ubv\]
and so
 $a\,\dee'\, ubv$ in $S$, giving that $a\,\leqj\, b$.
\end{proof}

It follows from Lemma ~\ref{jdonS} that if $S$ is a straight
left order in $Q$ inducing $(\leql,\leqr)$, then $\dee'$ is
an equivalence relation and $\leqj$ is a preorder with associated
equivalence relation $\jay'$.

\section{Partial orders of semigroups}\label{po}

Let $P$ be a partially ordered set and let $S$ be a semigroup.
We say that {\em $S$ is a partial order $P$ of subsets $S_{\alpha}$}
if $S$ is the disjoint union of
non-empty subsets $S_{\alpha}$, $\alpha\in P$ such
that for any $\alpha,\beta\in P$,
\[S_{\alpha}S_{\beta}\subseteq \bigcup\{ S_{\gamma}:\gamma\leq \alpha
\mbox{ and }\gamma\leq \beta\}.\]Equivalently,
there is a preorder $\preceq$ on $S$ satisfying the condition 
that for all $a,b\in S$, $ab\preceq a$
and $ab\preceq b$,  such that the classes
of the associated equivalence relation on $S$ are indexed by $P$.

Our definition is, of course, motivated by that of a semilattice
of semigroups. If $P$ is a semilattice and $S$ is a semilattice
$P$ of subsemigroups $S_{\alpha},\alpha\in P$, then certainly 
$S$ is a partial order $P$ of the subsets $S_{\alpha},\alpha\in P$.
However, there are certainly natural decompositions of semigroups
as partial orders of subsets that are not
necesarily associated with a semilattice
decomposition.

\begin{Example}\label{thejorder} The $\jay$-class
decomposition.
\end{Example}\begin{proof}
For any elements $a,b$ of a semigroup $S$ we have
$ab\, \leq_{\mathcal{J}}\, a$ and $ab\, \leq_{\mathcal{J}}\, b$.
\end{proof}

\begin{Example}\label{thejaysorder} The $\jay^*$-class decomposition.
\end{Example}
\begin{proof} We recall from \cite{fountainabundant} that
for an element $a$ of a semigroup $S$, there is a smallest
$*$-ideal $J^*(a)$
containing $a$, where an ideal $I$ of $S$ is a 
$*$-ideal if it is a union of $\ars$-classes and of 
$\els$-classes. 

The relation $\leqjays$ is defined on $S$ by the rule that
$a\,\leqjays\, b$ if and only if $J^*(a)\subseteq J^*(b)$.
It is easy to see that $\leqjays$ is a preorder, with associated
equivalence $\jays$.
It follows from (3) of
 Lemma 1.7 of \cite{fountainabundant},
that
for any $a,b\in S$ we have
$ab\,\leq_{\mathcal{J}^*}\, a$
and
$ab\,\leq_{\mathcal{J}^*}\, b$.
\end{proof}

We remark that in a regular semigroup, every ideal is a $*$-ideal.
Consequently, $\leqjay=\leqjays$ and $\jay=\jays$. Example
~\ref{thejorder} is thus a special case of Example~\ref{thejaysorder}.

Suppose now that $S$ is a semigroup and a partial order $P$ of
subsets $S_{\alpha},\alpha\in P$. For each $\alpha\in P$ 
we associate two ideals of $S$,
\[I^S_{\alpha}=\cup\{ S_{\beta}:\beta<\alpha\}\mbox{ and }
J^S_{\alpha}=\cup\{ S_{\beta}:\beta\leq\alpha\}\]
(we allow the possibility that $I^S_{\alpha}=\emptyset$).
We can thus form the Rees quotient
\[\overline{S_{\alpha}}=J^S_{\alpha}/I^S_{\alpha}=S_{\alpha}\cup
\{I^S_{\alpha}\},\]
or, in the case where $I^S_{\alpha}=\emptyset$,
we put $\overline{S_{\alpha}}=S_{\alpha}$.
We can think of  these Rees quotients
as 
`slices'
describing the local structure of $S$, 
and of $S$ as being built
up as `layers' of these slices.  In the next sections
we consider 
the relationship between 
$S$ being a straight left order
and the Rees quotients $\overline{S_{\alpha}},\alpha\in P$,
being straight left orders.

\section{Slicing}\label{slicing}

Throughout this section,   $S$ denotes a straight left order in a semigroup
$Q$, where $Q$ is a partial order $P$ of subsets
$Q_{\alpha},\alpha\in P$. For $\alpha\in P$ put
\[S_{\alpha}=Q_{\alpha}\cap S.\]
Notice that as $\overline{Q_{\alpha}}$ is a quotient of ideals
of $Q$, the set of non-zero elements $Q_{\alpha}$ is a union
of $\jay$-classes of $Q$. Since $S$ intersects every $\eh$-class
of $Q$, certainly each $S_{\alpha}\neq\emptyset$.
Clearly then $S$ is then a partial order $P$ of subsets $S_{\alpha},
\alpha\in P$.

The semigroup $Q$ is regular, consequently,
 for non-zero elements
$p,q$ of $\overline{Q_{\alpha}}$,
\[q\,\leqel \, p\mbox{ in }\overline{Q_{\alpha}}
\mbox{ if and only if }q\, \leqel\, p\mbox{ in }Q\]
and
\[q\,\leqar\, p\mbox{ in }\overline{Q_{\alpha}}
\mbox{ if and only if }q\,\leqar\, p\mbox{ in }Q.\]

Let $(\leql,\leqr)$ denote the $*$-pair
on $S$ induced by $Q$.
Since each $Q_{\alpha}$ is a union of $\jay$-classes,
we must have that $S_{\alpha}$ is a union of
$\el'$-classes and of $\ar'$-classes.

We remark that for each $\a\in P$, $\overline{S_{\alpha}}$ may be regarded
as a subsemigroup of $\overline{Q_{\alpha}}$. Clearly
 $I^Q_{\alpha}=\emptyset$
if and only if $I^S_{\alpha}=\emptyset$. In the case 
where $I^Q_{\alpha}\neq\emptyset$, we identify
the zeroes $\{ I^Q_{\alpha}\}$ and $\{ I^S_{\alpha}\}$
of $\overline{Q_{\alpha}}$ and $\overline{S_{\alpha}}$.

\begin{Lemma}\label{sliceup}
For each $\alpha\in P$,
 $\overline{S_{\alpha}}$ is a straight weak
left order
in $\overline{Q_{\alpha}}$.
\end{Lemma}
\begin{proof} 
If $q\in\overline{Q_{\alpha}}$, then either $q=0=0\s 0$,
or $q\in Q_{\alpha}$ and $q=a\s b$ where $a,b\in S$ and $a\,\ar\, b$ in
$Q$. Since $q\,\eh\, b$ in $Q$ we must have that $a,b\in S_{\alpha}$
and $a\,\ar\, b$ in $\overline{Q_{\alpha}}$.
\end{proof}

From Lemma~\ref{sliceup} we certainly have the following.

\begin{Corollary}\label{slicedorders} For any $\alpha\in P$,
$\overline{S_{\alpha}}$ is a straight left order in $\overline{Q_{\alpha}}$
if and only if for any $a\in S_{\alpha}$,
\[a\,\ehs\, a^2\mbox{ in }\overline{S_{\alpha}}\mbox{ implies that }
a\,\ehs\, a^2\mbox{ in }S.\]
\end{Corollary}

In general we do not
have enough control over the square cancellable elements,
for the condition in Corollary ~\ref{slicedorders} to hold.
However, in two cases of note we have the required result.

\begin{Proposition} \label{abundantslices}
Let $\a\in P$. If 
 $Q$ is completely regular
or $S$ is abundant and $J^S_{\alpha}$ and $I^S_{\alpha}$
are $*$-ideals,
then $\overline{S_{\alpha}}$ is a straight left order in
$\overline{Q_{\alpha}}$.
\end{Proposition}
\begin{proof} We need only consider (ii). In this case
the condition required in Corollary~\ref{slicedorders}
follows from Lemma 2.5 of \cite{basisii}.
\end{proof}

The canonical decomposition of a semigroup of left quotients $Q$
as a partial order $P$ of subsets $Q_{\alpha},\a\in P$, arises
from the preorder $\leqjay$, as in Example~\ref{thejorder}. In this
case the Rees quotients $\overline{Q_{\alpha}}$ are, of
course, the principal factors of $Q$, and we abuse terminology
by saying `$P$ is the $\leqjay$-order on $Q$'.
The canonical decomposition of a  left order $S$
as a partial order $P$ of subsets $S_{\alpha},\a\in P$, arises
from the preorder $\leqjays$, as in Example~\ref{thejaysorder}; in this
case we say `$P$ is the $\leqjays$-order on $S$'. If $P$
arises from the preorder
$\leqj$ we say `$P$ is the $\leqj$-order'. The relations
$\leqjay$ on $Q$ and $\leqjays$ and
$\leqj$ on $S$ are, of course, related,
and we see examples of this in the following proposition and
in Sections~\ref{fs} and ~\ref{cs}.

\begin{Proposition}
Let $S$ be a straight left order in $Q$ where $Q$ induces
$(\leql,\leqr)$. Suppose that $Q$ is a partial order
$P$ of subsets $Q_{\alpha},\a\in P$,
so that $S$ is a partial order $P$ of subsets
$S_{\alpha}=S\cap Q_{\alpha},\a\in P$. Then $P$ is
the $\leqjay$-order on $Q$ if and only if
$P$ is the $\leqj$-order on $S$.

Moreover, if $S$ is stratified in $Q$, that is
$\el'=\els$ and $\ar'=\ars$, then $P$ is the
$\jay$-order on $Q$ if and only if $P$ is the
$\jays$-order on $S$.

\end{Proposition}
\begin{proof}
The first part follows from Lemma~\ref{jdonS} and the facts
that $S$ intersects every $\eh$-class of $Q$
and each $Q_{\alpha}$ is a union of $\eh$-classes.

\medskip

For the remainder of the proof, we assume that $S$ is stratified
in $Q$. 

\medskip

Suppose that $P$ is the $\leqjay$-order on $Q$, so that
$P$ is the $\leqj$-order on $S$. Let $a,b\in S$ and
suppose that $a\,\leqj\, b$, so that
$a\,\dee'\, ubv$ for some $u,v\in S$. We have
that $\dees=\dee'$, so that by Lemma
1.7 of \cite{fountainabundant}, $a\,\leqjays\, b$ in $S$. 
On the other hand, if $c,d\in S$ and $c\,\leqjays\, d$, then
again by Lemma 1.7 of \cite{fountainabundant}, there
are elements $c_0,c_1,\hdots ,c_n\in S$,
say $c_i\in S_{\alpha_i}$, and
$x_1,\hdots ,x_n,y_1,\hdots ,y_n\in S^1$ such
that $d=c_0,c=c_n$ and
$(c_i,x_ic_{i-1}y_i)\in \dees$
for $i=1,\hdots ,n$. Since $\dee'=\dees$ and each
$S_{\alpha}$ is a union of
$\dees$-classes, we have
\[\a_0\geq\a_1\geq\hdots \geq \a_n.\]
Thus $c\,\leqjay\, d$ in $Q$ so that $c\,\leqj\, d$ in $S$.
Consequently,
 $P$ is the $\leqjays$-order on $S$.

Finally, suppose that $P$ is the $\leqjays$-order on $S$. 
Let $p\in Q_{\alpha}$ and $q\in Q_{\beta}$, and choose
$a\in S_{\alpha},b\in S_{\beta}$ with $p\,\eh\, a$ and $q\,\eh\, b$. 

If $\alpha\leq \beta$, then $a\,\leqjays\, b$, so that
by Lemma 1.7 of \cite{fountainabundant}, there are elements
$a_0,a_1,\hdots ,a_n\in S$, $x_1,\hdots ,x_n,y_1,\hdots ,y_n\in S^1$
such that
$b=a_0,a=a_n$ and
$(a_i,x_ia_{i-1}y_i)\in \dees$ for
$i=1,\hdots ,n$. Hence in $Q$, $a_i\,\dee\, x_ia_{i-1}y_i$
so that   $a_i\,\dee'\, x_ia_{i-1}y_i$ and
\[a=a_n\,\leqj\, a_{n-1}\,\leqj\, \hdots \leqj a_0=b.\]
It follows that $a\,\leqjay\, b$ in $Q$, whence
$p\,\leqjay\, q$ in $Q$.

Conversely, if $p\,\leqjay\, q$, then $a\,\leqj\, b$ in
$S$ so that $a\,\dee'\, ubv$ for some $u,v\in S$,
giving that $\alpha\leq\beta$ as required.
\end{proof}

\section{Layering}\label{layering}

Suppose that a semigroup $S$ is a partial order $P$ of
subsets $S_{\alpha},\a\in P$ such that
each $\overline{S_{\alpha}}$ is a straight left order
in a semigroup $W_{\alpha}$. The aim of this section is
to give necessary and sufficient conditions on $S$ such
that $S$ is a straight left order in $Q$ and $Q$
is a partial order $P$ of subsets $Q_{\alpha}$
such that $\overline{Q_{\alpha}}=W_{\alpha}$ for each $\a\in P$.

We  vary standard notation slightly by letting
$\leq_{\ell}^{\alpha}$ and $\leq_r^{\alpha}$ be the restriction to $S_{\alpha}$
of the preorders induced on $\overline{S_{\alpha}}$ by
$\leqel$ and $\leqar$ in 
$W_{\alpha}$ (that is, we exclude from the preorders
pairs of the form $(I^S_{\alpha},x)$).
We  denote by $\el'_{\alpha}$ and $\ar'_{\alpha}$ the equivalences
on $S_{\alpha}$ induced by $\leq_{\ell}^{\alpha}$ and $\leq_r^{\alpha}$
respectively; thus $\el'_{\alpha}$ and $\ar'_{\alpha}$
are the restrictions of Green's relations $\el$ and $\ar$
on $W_{\alpha}$ to the elements of $S_{\alpha}$.

\begin{Definition}\label{lr} The relations $\leql,\leqr$ are
defined on the
elements of $S$ by the rule that
 for any $a,b\in S$,
\[a\,\leql\, b\mbox{ if and only if }a, cb\in S_{\alpha}\mbox{
and }a\,\el_{\alpha}'\, cb\]
for some $\a\in P$. The relation $\leqr$ is defined dually.
\end{Definition}

In the above definition, if $a,cb\in S_{\alpha}$
and $a\,\el'_{\alpha}\, cb$, then choosing a square cancellable
element $x$ of $\overline{S_{\alpha}}$ in the $\ar'_{\alpha}$-class of $cb$,
we have $a\,\el'_{\alpha}\, xcb$ and $xc\in S_{\alpha}$ as $J^S_{\alpha}$
and $I^S_{\alpha}$ are ideals.
We may then immediately deduce
 the following.

\begin{Lemma}\label{erlemma} Let $a,b\in S_{\alpha}$. Then
\[a\,\leql\, b\mbox{ if and only if }a\,\leq_{\ell}^{\alpha}\, b\]
and dually,
\[a\,\leqr\, b\mbox{ if and only if }a\,\leq_{r}^{\alpha}\, b.\]
\end{Lemma}

\begin{Lemma} For any $a,b\in S$,
\[a\,\leql\, b\,\leql\, a,\]
if and only if   $a,b\in S_{\alpha}$ for
some $\a\in P$ and $a\,\el'_{\alpha}\, b$. The dual result holds for
$\leqr$.
\end{Lemma}
\begin{proof}
Suppose that $a,b\in S$ and
\[a\,\leql\, b\,\leql\, a.\]
We have that $a,cb\in S_{\alpha}$ and $b,da\in S_{\beta}$
for some $\a,\beta\in P$ and $c,d\in S$.
It follows from the definition of partial order of subsets
that $\a=\beta$ and
and so $a\el'_{\alpha} b$ by Lemma~\ref{erlemma}.

Conversely, if $a,b\in S_{\alpha}$ and $a\,\el'_{\alpha}\, b$,
then $a\leq_{\ell}^{\alpha} b$
and $b\leq_{\ell}^{\alpha} a$. By Lemma ~\ref{erlemma} we have
that
\[a\,\leql\, b\leql\, a\]
as required.
\end{proof}

Consequent upon the previous result, if $\leql$ and $\leqr$ are
transitive, then the associated equivalence relations
$\el'$ and $\ar'$ are
respectively the union of the equivalence relations $\el'_{\alpha}$
and $\ar'_{\alpha},
\a\in P$.

Recall that if $T$ is a left  order in semigroups $Q$ and $W$,
then $Q$ is {\em isomorphic to }$W$ {\em over }$S$ if there
is an isomorphism from $Q$ to $W$ that restricts to the identity
on $S$.

\begin{Theorem}\label{main}

 The following conditions on
the semigroup $S$ are equivalent:

(I) $(\leql,\leqr)$ is an embeddable $*$-pair and
$\mathcal{G}(S)=\sqs$;

(II) $(\leql,\leqr)$ is a $*$-pair such that

  for all $a,b,c\in S$ with $a\in \sqs$,
\[a^2b\,\ar'\, a^2c\mbox{ implies that }ab\,\ar'\, ac\]
and
$\sqs=\mathcal{G}(S)$;

(III) $S$ is a straight left order in a semigroup
$Q$ inducing
$(\leq_l,\leq_r)$.

If any of these conditions hold, then $Q$ is
 a partial order $P$ of subsets $Q_{\alpha},\a\in P$,
such that
each $\overline{S_{\alpha}}$ is a straight left order in 
$\overline{Q}_{\alpha}$
and $\overline{Q_{\alpha}}$ is
isomorphic to $W_{\alpha}$ over $\overline{S_{\alpha}}$.
\end{Theorem}
\begin{proof} {\em (III) implies (I)} This is immediate from
Theorem 4.4 of \cite{giiigiv}.

{\em (I) implies (II)} Suppose that
$a,b,c\in S$ with $a\in\sqs$ and $a^2b\,\ar'\, a^2c$. 
Since $(\leql,\leqr)$ is an embeddable $*$-pair, we
have $ab,ac\leqr a$ by condition (Eii)(r),
 so that by (Evii)(r),
$ab\,\ar'\, ac$.

{\em (II) implies (III)}
 By the construction of the relation $\eh'$
 condition G(ii) holds and we are given that
 G(i) also holds.

 We verify that $(\leql,\leqr)$ is an embeddable
$*$-pair. Since each $S_{\alpha}$ is
a union of $\el'$-classes and $\ar'$-classes and $\el'_{\alpha}\circ
\ar'_{\alpha}=
\ar'_{\alpha}\circ\el'_{\alpha}$ we have that (Ei) holds.
Conditions (Eii)(l) and (Eii)(r)
follow from Definition 3.1 and the remarks preceding the statement of
the theorem. Clearly
(Eiii) holds.

We proceed via a series of lemmas to verify the remaining
conditions.

\begin{Lemma} Let $a\in\mathcal{G}(S)$ and let $b\in S$ with
$b\,\leql\, a$. Then $b\,\ar'\, ba$.
\end{Lemma}
\begin{proof} We have that $b\,\el'_{\alpha} ca$ for some $\a\in P$
and $b,c\in S_{\alpha}$. Now $\ar'$ is a left congruence so that
$ca\,\ar'\, ca^2$ and in
$\overline{S_{\alpha}}$, 
\[b=qca\,\ar\, qca^2=ba\]
for some $q\in W_{\alpha}$. Hence $b\,\ar'_{\alpha}\, ba$ and so
$b\,\ar'\, ba$ as required.
\end{proof}

The above lemma show that condition (Ev)(l) holds; the proof
for (Ev)(r) is dual.

\begin{Lemma}\label{(Evii)(r)} Let $a\in\mathcal{G}(S)$ and let $b,c\in S$ with
$b,c\, \leqr\, a$ and $ab\,\ar'\, ac$. Then $b\,\ar'\, c$.
\end{Lemma}
\begin{proof}  By condition (Eii)(r) we have
that $b\,\ar'\, ax$ and $c\,\ar'\, ay$ for some 
$x,y\in S$. Since $\ar'$ is a left congruence,
\[a^2x\,\ar'\, ab\,\ar'\, ac\, \ar'\, a^2y\]
so that by the given hypotheses of (II) we have that
$ax\,\ar'\, ay$ and so $b\,\ar'\, c$.
\end{proof}

\begin{Lemma}\label{(Evi)(r)} Let $a\in\mathcal{G}(S)$ and let $b,c\in S$ with
$b,c\,\leqr\, a$ and $ab=ac$. Then $b=c$.
\end{Lemma}
\begin{proof} By condition (Ev)(r) we have
\[b\,\el'\, ab=ac\,\el'\, c.\]
Hence $b\,\el'_{\alpha}\, ab$ for some $\a\in P$. 
Since $\overline{S_{\alpha}}$ is a straight left order in $W_{\alpha}$
we may apply Lemma 4.7 of \cite{giiigiv} to deduce the
existence of $h\in \sqs,k\in S$ with $b\,\ar'\, h\,\ar'\, k$ and 
\[hb=kab.\]
We know that $\ar'\subseteq\ars$ so that $hc=kac$ and
\[hb=kab=kac=hc.\]
By Lemma~\ref{(Evii)(r)} we have that $b\,\ar'\, c$ and so
in $W_{\alpha}$,
\[b=h\s hb=h\s hc=c.\]
\end{proof}

\begin{Lemma}\label{(Evii)(l)} Let $a\in\mathcal{G}(S)$ and $b,c\in S$
with $b,c\,\leql\, a$ and $ba\,\el'\, ca$. Then $b\,\el'\, c$.
\end{Lemma}\begin{proof} By condition (Eii)(l)
we have that 
\[b\,\el'\, xa\mbox{ and }c\,\el'\, ya\]
for some $x,y\in S$. Since $\el'$ is a right congruence
\[xa^2\,\el'\, ba\,\el'\, ca\,\el'\, ya^2\]
and so 
\[xa^2\,\el'_{\alpha}\, ya^2\]
for some $\a\in P$. 
By Lemma 4.7 of \cite{giiigiv}, 
\[hxa^2=kya^2\]
for some $h\in\sqs,k\in S$ and
\[h\,\ar'\, k\,\ar'\,  xa^2\,\ar'\, xa.\]
We know that $\ar'\subseteq \ars$ and so 
\[hxa=kya,\]
giving by condition (Ev)(r) 
\[xa\,\el'\, hxa=kya\leql ya.\]
Hence $b\,\leql\, c$; dually, $c\,\leql\, b$ so that
$b\,\el'\, c$ as required.
\end{proof}

\begin{Lemma}\label{(Evi)(l)} Let $a\in\mathcal{G}(S)$ and let $b,c\in S$ with
$b,c\,\leql\, a$ and $ba=ca$. Then $b=c$.
\end{Lemma}
\begin{proof} From Lemma ~\ref{(Evii)(l)} we have that
$b\,\el'\, c$, say $b\,\el'_{\alpha}\, c$ and so
by Lemma 4.7 of \cite{giiigiv}, $hb=kc$ for some
$h\in\sqs,k\in S$ with $h\,\ar'\, k\,\ar'\, b$. By condition
(Ev)(l) 
\[b\,\ar'\, ba=ca\,\ar' c.\]
From
\[hba=kca=kba\]
we deduce that
$hb=kb$ and $hc=kc$. Hence
\[hb=kc=hc.\]
In $W_{\a}$ we then have
\[b=h\s hb=h\s hc=c.\]
\end{proof}

Lemmas ~\ref{(Evi)(l)},\ref{(Evi)(r)},\ref{(Evii)(l)} and
~\ref{(Evii)(r)} show that conditions (Evi)(l), (Evi)(r),
(Evii)(l) and (Evii)(r) hold, respectively.

We call upon Theorem 4.4 of \cite{giiigiv} 
to deduce that $S$ is  a straight left order in $Q$ where
$Q$ induces $(\leql,\leqr)$.

\medskip

Suppose now that any of the equivalent conditions (I) to (III) holds.
We show that $Q$ is a partial order $P$ of subsets $Q_{\alpha}$, $\a\in P$.

Let $q\in Q$. Since $S$ is a straight left order in $Q$,
$q=h\s k$ where $h\in\sqs,k\in S$ and $h\,\ar'\, k$ in $S$.
Hence $h,k\in S_{\alpha}$ for some $\a\in P$ and $q\,\eh\, k$ in $Q$.
Thus
\[Q=\bigcup\{ Q_{\alpha}:\a\in P\}\]
where
\[Q_{\alpha}=\{ q\in Q:q\,\eh'\, k,k\in S_{\alpha}\}.\]
Since each $S_{\alpha}$ is a union of $\eh'$-classes,
$Q_{\alpha}\cap Q_{\beta}=\emptyset$ for $\alpha\neq \beta$
and $S_{\alpha}=S\cap Q_{\alpha}$.

Let $p\in Q_{\alpha},q\in Q_{\beta }$. Then $p\,\eh\, h$,
$q\,\eh\, k$ for some $h\in S_{\alpha},k\in S_{\beta}$.
Now
\[pq\,\ar\, pk\,\el\, hk\]
and $hk\in S_{\gamma}$ for some $\gamma\leq \a,\beta$.
Since each $S_{\gamma}$ is a union of $\ar'$-classes and
of $\el'$-classes, it follows that $pq\in Q_{\gamma}$.
Hence $Q$ is a partial order of the subsets $Q_{\alpha},\a\in P$.

For each $\a\in P$ we let $\overline{Q_{\alpha}}$ be the Rees quotient
associated with $Q_{\alpha}$. From Corollary ~\ref{slicedorders}
 and the construction
of the relation $\eh'$, we have that $\overline{S_{\alpha}}$ is a straight
left order in $\overline{Q_{\alpha}}$.

Let $a,b\in \overline{S_{\alpha}}$. If $a=\{ I^S_{\alpha}\}$
then $a$ is the zero of both
$W_{\alpha}$ and $\overline{Q_{\alpha}}$,
so that $a\leqar b$ in both $W_{\alpha}$ and
$\overline{Q_{\alpha}}$.
Otherwise, $a\in S_{\alpha}$ and
\[\begin{array}{rcl}
a\,\leqar\, b\mbox{ in }W_{\alpha}&\Leftrightarrow &b\in S_{\alpha}\mbox{ and }
a\,\ar'\, bc,\mbox{ for some }c\in S_{\alpha}\\
&\Leftrightarrow& 
b\in S_{\alpha}\mbox{ and }
a\,\ar'\, bc,\mbox{ for some }c\in S\\
&\Leftrightarrow& b\in S_{\alpha}\mbox{ and }a\leqr b\\
&\Leftrightarrow& b\in S_{\alpha}\mbox{ and }a\leqar b\mbox{ in }Q\\
&\Leftrightarrow& b\in S_{\alpha}\mbox{ and }a\leqar b\mbox{ in }
\overline{Q_{\alpha}}.\end{array}\]
Corollary 4.3 of \cite{uniqueness} gives that $\overline{Q_{\alpha}}$ is
isomorphic to $W_{\alpha}$ over $\overline{S_{\alpha}}$.
\end{proof}

The rest of the paper is devoted to specialisations of Theorem
~\ref{main}.

\section{Semilattices of straight left orders}\label{semilattices}
 Let $S$ be a semilattice $Y$ of subsemigroups $S_{\alpha}$, $\a\in Y$,
where each $S_{\alpha}$ is a straight left order in $T_{\alpha}$
inducing $(\leq_{\ell}^{\alpha},\leq_r^{\alpha})$ with
associated equivalence relations $\mathcal{L}'_{\alpha}$
and $\mathcal{R}'_{\alpha}$. As
remarked in Section~\ref{po}, $S$ is
certainly a partial order $Y$ of the subsets
$S_{\alpha},\a\in Y$.
We define the relations $\leql$ and $\leqr$ as in Definition ~\ref{lr}.
Let $\a\in P$. If $I^S_{\alpha}=\emptyset$, then we
let $W_{\alpha}=T_{\alpha}$.
If $I^S_{\alpha}\neq\emptyset$ put $W_{\alpha}=T_{\alpha}^0$.
In view of Lemma~\ref{zeroesagain}, in either case
$\overline{S_{\alpha}}$
is a straight left order in $W_{\alpha}$. It is easy to see that the pair
$(\leql,\leqr)$ arising from
the left quotient semigroups
$T_{\alpha}$  coincides with the
pair of relations arising  from
the left quotient semigroups $W_{\alpha}$. We can thus
apply Theorem~\ref{main} to obtain the first part of the following
result.

\begin{Theorem}\label{sl} The following conditions on the
semigroup $S$ are equivalent:

(A) $(\leql,\leqr)$ is a $*$-pair such that
for all $a\in\sqs,b,c\in S$,
\[a^2b\,\ar'\, a^2c\mbox{ implies that }ab\,\ar'\, ac\]
and $\sqs=\mathcal{G}(S)$;

(B) $S$ is a straight left order in $Q$ inducing
$(\leql,\leqr)$.

If any of these conditions hold, then $Q$ is a semilattice
$Y$ of subsemigroups
$Q_{\alpha},\alpha\in Y$, such that
each $S_{\alpha}$ is a straight left order in $Q_{\alpha}$ and
$Q_{\alpha}$ is isomorphic to $T_{\alpha}$ over $S_{\alpha}$.
\end{Theorem}
\begin{proof} By Theorem~\ref{main},
(A) and (B) are equivalent.
  If (A) or (B) hold, then again by Theorem ~\ref{main},
$Q$ is a partial order $Y$ of subsets $Q_{\alpha}$, $\a\in Y$,
such that each $\overline{S_{\alpha}}$ is a straight left
order in $\overline{Q_{\alpha}}$ and
$\overline{Q_{\alpha}}$ is isomorphic to $W_{\alpha}$ over
 $\overline{S_{\alpha}}$.

Let $\a\in P$. If $I^S_{\alpha}=\emptyset$, then as in
Section~\ref{slicing}, $I^Q_{\alpha}=\emptyset$. In this
case, $\overline{S_{\alpha}}=S_{\alpha}$,
 $\overline{Q_{\alpha}}=Q_{\alpha}$
is a subsemigroup of $Q$, and 
$W_{\alpha}=T_{\alpha}$. Hence
$S_{\alpha}$ is a straight left order
in $Q_{\alpha}$ and
$Q_{\alpha}$ is isomorphic to $T_{\alpha}$ over
$S{\alpha}$.

Otherwise, $\overline{S_{\alpha}}=S^0_{\alpha}$. By
Lemma~\ref{zeroes}, the zero of $\overline{S_{\alpha}}$ is
the zero of $\overline{Q_{\alpha}}$, and as
$\overline{S_{\alpha}}$ has no zero divisors, neither
does $\overline{Q_{\alpha}}$. Thus $Q_{\alpha}$ is a
subsemigroup of $Q$ and by Lemma~\ref{zeroesagain}, 
$S_{\alpha}$ is a straight left order in $Q_{\alpha}$. Now
$\overline{Q_{\alpha}}$ is isomorphic to $T^0_{\alpha}$
over $S^0_{\alpha}$; clearly then
$Q_{\alpha}$ is isomorphic to $T_{\alpha}$ over $S_{\alpha}$.

Finally, let $p\in Q_{\alpha}$ and $Q\in Q_{\beta}$. Then
$p\,\eh\, a$ and 
$q\,\eh\, b$ for some $a\in S_{\alpha}$
and $b\in S_{\beta}$, whence
\[pq\,\ar\, pb\,\el\, ab\]
so that as each $Q_{\gamma}$ is a union of $\dee$-classes, we have
$pq\in Q_{\alpha\beta}$. Thus $Q$ is
a semilattice $Y$ of subsemigroups $Q_{\alpha},\a\in Y$.
\end{proof}

Theorem ~\ref{sl}  has a number of applications. Theorem 5.5
of \cite{compreg} (see also \cite{paula})
 is the specialisation of Theorem~\ref{sl}
to the case where each $S_{\alpha}$ is a straight left order
in a completely simple semigroup $T_{\alpha}$. Surprisingly,
the conditions in (A) are only slightly stronger
than those in Theorem 5.5 of \cite{compreg}. In the latter
result one can use the decomposition of each $S_{\alpha}$ as a matrix
of right reversible cancellative subsemigroups to deduce that
$(\leql,\leqr)$
is a $*$-pair from the assumptions that $\leql$ and $\leqr$
are transistive, $\ar'$ is a left
congruence, $\el'\subseteq \els$ and $\ar'\subseteq \ars$
\cite{compreg}.
Moreover for any $a\in S$, $a\,\eh'\, a^2$ so that as $\eh'\subseteq
\ehs$,
$a\,\ehs\, a^2$. thus $\sqs=\mathcal{G}(S)$
holds automatically. Of course one can specialise further
to the case where each $T_{\alpha}$ is a group, to obtain the early
result of \cite{clifford} mentioned in the Introdution.

Specialisations of Theorem~\ref{sl} to the case where $S$ is
cancellative or commutative can be found in
\cite{cancellative} and \cite{commutative}.

\section{The fully stratified case}\label{fs} 

Recall that a straight left order $S$ in $Q$ is
{\em fully stratified} in $Q$ if $Q$ induces
$(\leqels,\leqars)$. Clearly we may specialise
Theorem~\ref{main} to fully stratified straight left
orders. The cleanest result is available when $S$ is IC
abundant.

\begin{Theorem}\label{fsic} Let $S$ be an idempotent connected abundant
semigroup. 

If $S$ is a fully stratified straight left order in $Q$,
where $Q$ is a partial order $P$ of subsets $Q_{\alpha},\a\in P$,
then for any $\a\in P$, $\overline{S_{\alpha}}$ is
an abundant fully stratified straight left order in
$\overline{Q_{\alpha}}$.

Conversely, suppose that $S$ is a partial order $P$ of
subsets $S_{\alpha}\a\in P$, such that each $\overline{S_{\alpha}}$
is an abundant fully stratified straight left order.
Then $S$ is a fully stratified straight left order if and only if
for all $a\in \sqs$ and $b,c\in S$,
\[a^2b\,\ars\, a^2c\mbox{ implies that }ab\,\ars\, ac.\]
\end{Theorem}
\begin{proof} Suppose that $S$ is a fully stratified
straight left order in $Q$, where $Q$ is a partial order
$P$ of subsets $Q_{\alpha},\a\in P$. As in Section~\ref{slicing},
$S$ is a partial order $P$ of subsets $S_{\alpha}=S\cap Q_{\alpha}$,
and each $S_{\alpha}$ is a union of $\els$-classes and of
$\ars$-classes. Thus for any $\a\in P$, $J^S_{\alpha}$ and
$I^S_{\alpha}$ are
$*$-ideals. From Proposition ~\ref{abundantslices}, each $\overline{S_{\alpha}}$
is a straight left order in $\overline{Q_{\alpha}}$. It remains to show
that $\overline{S_{\alpha}}$ is abundant and fully stratified
in $\overline{Q_{\alpha}}$. 

Let $a,b\in \overline{S_{\alpha}}$ and suppose that $a\,\leqars\, b$
in $\overline{S_{\alpha}}$. Either $a=0$,
in which
case certainly $a\,\leqar\, b$ in $\overline{Q_{\alpha}}$.
Otherwise, $a,b\in S_{\alpha}$ and by Lemma 2.5 of
\cite{basisii}, $a\,\leqars\, b$ in $S$. But $S$
is fully stratified in $Q$, so that $a\,\leqar\, b$ in $Q$ and
hence $a\,\leqar \, b$ in $\overline{Q_{\alpha}}$. Together
with the dual argument for $\leqels$, we have that
$\overline{Q_{\alpha}}$ induces
$(\leqels,\leqars)$ on $\overline{S_{\alpha}}$.

Again from Lemma 2.5 of \cite{basisii}, the non-zero $\els$-classes
and $\ars$-classes of $\overline{S_{\alpha}}$ are  $\els$-classes
and $\ars$-classes of $S$ and thus contain idempotents since
$S$ is abundant. Thus $\overline{S_{\alpha}}$ is abundant.

Conversely, suppose that $S$ is a partial order
$P$ of subsets $S_{\alpha},\a\in P$, such that
each $\overline{S_{\alpha}}$ is an abundant fully stratified
straight left order in a semigroup $W_{\alpha}$. Let $a\in S_{\alpha},
b\in S_{\beta}$ and suppose that $a\,\ars\, b$ in $S$.
Since $\overline{S_{\alpha}}$ and $\overline{S_{\beta}}$
are abundant, there are idempotents $e\in \overline{S_{\alpha}}$
and $f\in \overline{S_{\beta}}$ with $a\,\ars\, e$ in
$\overline{S_{\alpha}}$
and $b\,\ars\, f$ in
$\overline{S_{\beta}}$.
Hence $ea=a$ and $fb=b$.
It follows that
$e\in S_{\alpha}$ and $f\in S_{\beta}$. Moreover from
$ea=a$ and $a\,\ars\, b$ in $S$ we have
$eb=b$, so that $\beta\leq\alpha$; dually, $\a\leq\beta$.
Hence each $S_{\alpha}$ is a union of
$\ars$-classes, and dually, of $\els$-classes. Consequently,
for any $\a\in P$, $I^S_{\alpha}$
and $J^S_{\alpha}$ are $*$-ideals.

Let $\leql$ and $\leqr$ be as in Definition~\ref{lr}. 

Let $a,b\in S$ and suppose that $a\,\leql\, b$. Then
for some $\a\in P$ and $c\in S$,
$a,cb\in S_{\alpha}$ and $a\,\el'_{\alpha}\, cb$. As
$\overline{S_{\alpha}}$ is fully stratified
we have that $a\,\els\, cb$ in
$\overline{S_{\alpha}}$ so that by Lemma 2.5
of \cite{basisii}, $a\,\els\, cb$ in $S$, so that
$a\,\leqels\, b$ in $S$.

Conversely, if $a\,\leqels\, b$ in $S$, then by Lemma~\ref{ic},
$a\,\els\, xb$ for some $x\in S$. Since each $S_{\alpha}$ is a union of
$\els$-classes, we have $a,xb\in S_{\alpha}$ for some $\a\in P$ and by
Lemma 2.5 of \cite{basisii}, $a\,\els\, xb$ in
$\overline{S_{\alpha}}$. Thus
$a\,\el'_{\alpha}\, xb$ so that
$a\,\leql\, b$. Hence
$\leql=\leqels$ and dually, $\leqr=\leqars$.
Clearly $(\leql,\leqr)$ is a $*$-pair
and $\sqs=\mathcal{G}(S)$. We can now apply
Theorem ~\ref{main} to deduce that $S$ is  a straight
left order in a semigroup $Q$ inducing $(\leqels,\leqars)$, that is,
$S$ is a fully stratified straight left
order, if and only if for all $a\in \sqs$ and $b,c\in S$,
\[a^2b\,\ars\, a^2c\mbox{ implies that }ab\,\ars\, ac.\]
\end{proof}

We now present an appication of this theorem. Theorem 5.12 of
\cite{basisii} states that the semigroup $\en_f A$ of
endomorphisms of finite rank of a stable basis algebra is a fully
stratified straight left order. The proof proceeds
by verifying all the conditions necessary for
$(\leqels, \leqars)$ to be an
embeddable $*$-pair, together with condition (Gii). Theorem
~\ref{fsic} gives an alternative approach. For further details
of basis algebras, we refer the reader to \cite{basisi,basisii}.

It is clear that
\[\en_f A=\bigcup\{ U_n:0\leq n\leq\rk A,n\in\mathbb{N}\}\]
where
\[U_n=\{ \a\in \en A: \rk A=n\}.\]
For any $\alpha,\beta\in \en A$, Lemma 5.1 of
\cite{basisii} gives that
\[\rk \alpha\beta\leq\mbox{ min }\{ \rk \alpha,\rk\beta\} \]
so that in our terminology, $\en_f A$ is a partial order
(indeed, a chain) $N$
of subsets $U_n,n\in N$, where
$N$ is the set of finite ranks of endomorphisms of $A$. Our Rees
quotients $\overline{U_{n}}$ are, in the terminology of
\cite{basisii}, the semigroups $S_n/S_{n-1}$. The latter are shown
in the final section of \cite{basisii} to be primitive abundant
semigroups in which all the non-zero idempotents are $\dee$-related.
Moreover the non-null $\ehs$-classes are right reversible.
It follows from \cite{fountainabundant} that $S_n/S_{n-1}$ is
isomorphic to a Rees matrix semigroup
$\mathcal{M}^0(H;I,\Lambda;P)$,
where $H$ is a right reversible cancellative semigroup
and $I$ and $\Lambda$ index the $\ars$-classes and $\els$-classes
of $S$, respectively.
It is easy to see that $\mathcal{M}^0(H;I,\Lambda;P)$
is a fully stratified straight left order
in $\mathcal{M}^0(G;I,\Lambda;P)$,
where $G$ is the group of left quotients of $H$.

From Theorem 3.9 and Corollary 5.2 of \cite{basisii},
$\en_f A$ is an IC abundant semigroup. We thus have
exactly the right circumstances to
apply Theorem~\ref{fsic}. Lemma 5.11 of \cite{basisii}
gives that for all $a\in \sqs$ and $b,c\in S$,
\[a^2b\,\ars\, a^2c\mbox{ implies that }ab\,\ars\, ac.\]
Consequently, $\en_f A$ is a fully stratified straight left
order.

\section{Straight left orders in completely semisimple
semigroups}\label{cs}

We recall from \cite{cp} that a semigroup is
{\em completely semisimple} if all its principal factors
are completely (0)-simple.
We end the paper with a new characterisation of straight left
orders in completely semisimple semigroups.

Let $S$ be a semigroup possessing a $*$-pair $(\leql,\leqr)$.
We say that $S$ has $M_L^*$ ($M_R^*$) if the
$\el'$-class ($\ar'$-classes) {\em within any
$\dee$-class} satisfy the descending chain condition.

\begin{Theorem} The following are equivalent for a semigroup $S$:

(i) $S$ is a straight left order in a completely semisimple
semigroup $Q$ inducing $(\leql,\leqr)$;

(ii) $S$ is a straight left order in a semigroup $Q$ inducing
$(\leql,\leqr)$ such that $\jay'=\dee'$ and
$M^*_L$ and $M_R^*$ hold;

(iii) $S$ is a partial order $P$ of subsets $S_{\alpha},\a\in P$, 
each $\overline{S_{\alpha}}$ is a (straight)
left order in a completely 0-simple 
(or completely simple) semigroup such that
 $(\leql,\leqr)$ is a $*$-pair,
  for all $a,b,c\in S$ with $a\in \sqs$,
\[a^2b\,\ar'\, a^2c\mbox{ implies that }ab\,\ar'\, ac\]
and
$\sqs=\mathcal{G}(S)$.
\begin{proof} {\em (i) implies (ii)} In $Q$ we have that
$\jay=\dee$ and so  by Lemma~\ref{jdonS}, $\jay'=\dee'$
on $S$. Suppose that $a_1,a_2,\hdots $ are $\dee'$-related
elements of $S$ such that
\[a_1\geq_{\ell}a_2\geq_{\ell}\hdots \]
Then in $Q$ $a_1,a_2,\hdots$ are $\dee$-related
and we have
\[a_1\geq_{\mathcal{L}}a_2\geq_{\mathcal{L}}\hdots \]
so that the same is true in the principal factors.
 But these are
completely (0)-simple, and so
\[a_1\,\el\, a_2\,\el\hdots\]
in $\overline{Q_{\alpha}}$ and hence in $Q$. This
gives 
\[a_1\,\el'\, a_2\,\el'\hdots .\]
Certainly then $M^*_L$ holds. Dually condition
$M_R^*$ holds.

{\em (ii) implies (i)} We have $\jay=\dee$ on $Q$ since
$\jay'=\dee'$ on $S$. 
Consider a principal factor $W$ of $Q$. 
Let $q_1,q_2,\hdots $ be non-zero elements
of $W$ such that
\[q_1\geq_{\mathcal{L}}q_2\geq_{\mathcal{L}}\hdots \]
in $W$. For each $i$ choose $a_i\in S$
such that $a_i\,\eh\, q_i$ in $Q$ and hence
in $W$. Now in $W$ and $Q$
\[a_1\geq_{\mathcal{L}}a_2\geq_{\mathcal{L}}\hdots \]
so that
\[a_1\geq_{\ell}a_2\geq_{\ell}\hdots\]
But the $q_i$'s are all $\jay$-related in $Q$, so that
the $a_i$'s are all $\dee'$-related in $S$. Since
$S$ has $M_L^*$, we have that for some $n$,
\[a_n\el'\, a_{n+1}\,\el'\, \hdots \]
whence
\[q_n\el\, q_{n+1}\,\el\, \hdots \]
in $Q$ and $W$. Hence $W$
has the descending chain condition on principal left ideals,
and dually for principal right ideals. By Theorem 3.3.2 of \cite{howie},
$W$ is completely (0)-simple.

{\em (i) implies (iii)} We have that $Q$ is a partial order
$P$ of subsets $Q_{\alpha},\a\in P$, where $P$ is the
$\leqjay$-order on $Q$, and $S$ is thus a partial order $P$ of subsets
$S_{\alpha}=S\cap Q_{\alpha}$, where $P$ is
the $\leqj$-order on $S$. By Lemma 2.1, each $\overline{S_{\alpha}}$ is
a straight
weak left order in $\overline{Q_{\alpha}}$. But if $a\in
\overline{S_{\alpha}}$
is square cancellable in $\overline{S_{\alpha}}$, then either $a=0$, so
lies in a subgroup of
$\overline{Q_{\alpha}}$, or as $a\,\ehs\, a^2$ in
$\overline{S_{\alpha}}$, $a^2\neq 0$. In the latter case it is well
known that $a$ lies in a subgroup of $\overline{Q_{\alpha}}$. Thus
$\overline{S_{\alpha}}$ is a straight left order in a completely
(0)-simple semigroup.

It is easy to verify that $(\leql,\leqr)$ coincides with
the pair of preorders constructed as in Section~\ref{layering} from the
straight left orders $\overline{S_{\alpha}}$. Clearly then (iii) holds. 

{\em (iii) implies (i)} This is immediate from Theorem~\ref{main}.
\end{proof}
\end{Theorem}

\end{document}